\documentclass[11pt]{article}   
\ifx\pdfoutput\undefined
 \usepackage[latin1]{inputenc}
 \usepackage[T1]{fontenc}
 \usepackage{url}
 \usepackage[dvips]{graphicx}
  
\else
 \pdfoutput=1\pdfcompresslevel=9\pdfadjustspacing=1
 \usepackage[pdftex,bookmarks,a4paper,colorlinks]{hyperref}
 \usepackage{aeguill}
 \usepackage[pdftex]{graphicx}
 \hypersetup{bookmarksnumbered,plainpages=false,hypertexnames=false}
 \hypersetup{pagecolor=blue,linkcolor=blue,citecolor=blue,urlcolor=blue}
 \hypersetup{pdfcreator=PDFLaTeX with TPAgregTeam class}
 \usepackage{thumbpdf}
\fi

\usepackage[english]{babel}
\usepackage{amsmath,amsfonts,amssymb}
\usepackage{theorem}
\usepackage{geometry}
\usepackage{enumerate}

\theoremstyle{plain}
\newtheorem{thrm}{Theorem}[section]
\newtheorem{lmm}[thrm]{Lemma}
\newtheorem{crllr}[thrm]{Corollary}

\theoremstyle{definition}
\newtheorem{dfntn}[thrm]{Definition}

\newtheorem{rmrk}[thrm]{Remark}

\newcommand{\proofend}{\hfill $\Box{~}$}
\newenvironment{proofp}[1]
               {\noindent {\emph{Proof {#1}.}}}
               {\proofend\\}

\newenvironment{proof}
               {\noindent {\textbf{Proof.}}}
               {\proofend\\}

\def\R{\mathbb R}

\def\dE{\mathbb E}
\def\dR{\mathbb R}
\newcommand{\cB}{\mathcal{B}}
\newcommand{\cF}{\mathcal{F}}
\newcommand{\cI}{\mathcal{I}}
\newcommand{\cL}{\mathcal{L}}
\newcommand{\cN}{\mathcal{N}}
\newcommand{\ABS}[1]{{\left| #1 \right|}} 
\newcommand{\PAR}[1]{{\left(#1\right)}} 
\newcommand{\SBRA}[1]{{\left[#1\right]}} 
\newcommand{\BRA}[1]{{\left\{#1\right\}}} 
\newcommand{\gd}{\Gamma{\!\!_2}}
\newcommand{\ep}[1]{\text{Ent}^\Phi_#1}
\newcommand{\enp}[2]{\text{Ent}^\Phi_#1\PAR{#2}}

\begin{document}

\begin{center}
  \begin{center}
    \begin{Large}
\textbf{\textsc{Logarithmic Sobolev Inequalities for Inhomogeneous Markov
    Semigroups}}    
    \end{Large}
  \end{center}
\end{center}
\bigskip
\bigskip
\begin{center}
  \begin{large}
\textsc{Jean--François Collet}\footnote{Laboratoire J.A. Dieudonné, Université de Nice Sophia
  Antipolis, Parc Valrose 06108 Nice Cédex 02}
And 
\textsc{Florent Malrieu}\footnote{IRMAR, Université
    Rennes 1, Campus de Beaulieu,
    35042 Rennes Cedex}    
  \end{large}
\end{center}
\bigskip
\bigskip

\date{\today}
\noindent
\textbf{Abstract.}  We investigate the dissipativity properties of a
class of scalar second order parabolic partial differential equations
with time-dependent coefficients. We provide explicit condition on the
drift term which ensure that the relative entropy of one particular
orbit with respect to some other one decreases to zero. The decay rate
is obtained explicitly by the use of a Sobolev logarithmic inequality
for the associated semigroup, which is derived by an adaptation of
Bakry's $\Gamma-$ calculus.  As a byproduct, the systematic method for
constructing entropies which we propose here also yields the
well-known intermediate asymptotics for the heat equation in a very
quick way, and without having to rescale the original equation.

\bigskip
\noindent
\textbf{Résum\'e.} Cet article propose une \'etude du m\'ecanisme de
dissipation d'entropie pour une classe d'équations aux dérivées
partielles paraboliques dont les coefficients dépendent du temps.
Sous des critères formul\'es explicitement en terme des coefficients,
nous établissons la décroissance exponentielle de l'entropie relative
d'une orbite par rapport à une autre, pour des \'equations n'admettant
pas de solution stationnaire. La m\'ethode utilis\'ee repose sur
l'obtention d'une inégalité de type Sobolev logarithmique pour le
semi-groupe associé, grâce à une adaptation du critère de Bakry-Émery.



\tableofcontents

\section{Introduction}\label{s1}

\subsection{The Kullback-Leibler Distance as a Particular $\Phi$-Entropy}

Given two probability densities $u,v$ on $\R^d$, the entropy of $u$
relative to $v$ (also known in information theory as their
Kullback-Leibler distance (see \cite{KL}), although it is not a
distance) is defined by
\begin{equation}
H(u|v) := \int_{\R^d} u(x) \ln \PAR{\frac{u(x)}{v(x)}} \, dx.
\label{log-ent}
\end{equation}
Although this quantity does not satisfy the triangle inequality, it is
always nonnegative and vanishes only when $u=v$. These two facts are
immediate consequences of the well-known Pinsker inequality (see
\cite{Pi}):
$$
H(u|v) \geq \frac{1}{2} |u-v|_{L^1}^2.
$$
Thus the quantity $H(u|v)$ may provide some notion of ``distance''
between $u$ and $v$. In Partial Differential Equations $H(u|v)$ may be
useful in studying the asymptotic behavior of a dissipative system. In
this context $v$ is in general a stationary solution (or in physical
terminology a detailed balance equilibrium), $u$ is the orbit of some
Kolmogorov (or some other parabolic) equation, and $H(u|v)$ is a
decreasing function of time. In many cases this fact can be combined
with some clever inequalities to show that $v$ in fact attracts $u$ in
some appropriate metric. This is the basis of the well-known entropy
dissipation method (see \cite{T2}), which has been used to good
advantage in many examples such as (linear or nonlinear) parabolic
equations, kinetic equations, etc...

As is well-known from information theory and Statistical Physics
(\ref{log-ent}) is a particular instance of a more general class of
entropies:
\begin{equation}
H^\Phi(u|v) := \int_{\R^d} v(x) \Phi \PAR{\frac{u(x)}{v(x)}} \, dx,
\label{phi-ent}
\end{equation}
where $\Phi$ is any convex function defined on $[0,\infty[$. Formula
(\ref{log-ent}) corresponds to the particular choice $\Phi(z)=z \log z$ which
has an interesting extensivity property \cite{C1}, but as a general rule the
dissipation of entropy is a convexity property which has little to do with the
specific properties of the $z \log z$ function.

In this paper we investigate a class of linear parabolic equations which due
to the presence of time-dependent coefficients have no stationary solution,
but for which orbits still do come together in the ``metric'' given by
\eqref{log-ent} or \eqref{phi-ent}. More precisely, we will give explicit
conditions on the coefficients which ensure that the quantities
\eqref{log-ent} or \eqref{phi-ent} decrease to zero for large time, with
quantitative bounds, for any two orbits $u$ and $v$, that is, {\bf even when
  $v$ is a non-stationary solution}.

\subsection{The Entropy Production for Linear Scalar Advection-Diffusion 
  Equations}

Let us consider a general multidimensional linear equation in the form:
\begin{equation}
\frac{\partial u}{\partial t} + \hbox{\rm div} J(u,\nabla u) = 0,
\label{adv-diff1}
\end{equation}
where $J$ is some specified flux function.  Let $\Phi$ be any convex
function, and assume that $J$ is such that this equation preserves
positivity. Given two (time-dependent or not) positive solutions $u$
and $v$ of this equation, we define the quantity $H^\Phi(u|v)$ at any time
$t$ by
$$
H^\Phi(u|v) := 
\int_{\R^d} v(t,x) \Phi \PAR{\frac{u(t,x)}{v(t,x)}} \,dx
$$
(From here on, the notation $dx$ will mean $d-$dimensional Lebesgue
measure).  The convexity of $\Phi$ implies a lower bound for $H^\Phi(u|v)$ as
follows: using Jensen's inequality with the probability measure
$\frac{vdx}{\int v \, dx}$ we obtain the inequality
$$
\int_{\R^d} v(t,x) \Phi \PAR{\frac{u(t,x)}{v(t,x)}} \, dx \geq 
\int_{\R^d} v(t,x) \, dx \Phi \PAR{\frac{\int_{\R^d} u(t,x) \, dx}{ 
\int_{\R^d} v(t,x) \, dx}}.
$$
If we assume that $\Phi(1)=0$ and that $u$ and $v$ have the same
integral (which we only have to assume at time zero, since equation
(\ref{adv-diff1}) will preserve mass for reasonable solutions), we see
that $H^\Phi(u|v)$ remains nonnegative.

Assuming that $u$ and $v$ vanish at infinity, a straightforward
computation then yields:
\begin{equation}
\frac{d}{dt} H^\Phi(u|v) = 
\int_{\R^d} \Phi''\PAR{\frac{u}{v}} \nabla \PAR{\frac{u}{v}} 
\SBRA{\frac{J(u, \nabla u)}{u} - \frac{J(v, \nabla v)}{v}} u \, dx.
\label{ent-prod}
\end{equation}
Perhaps the best known case of this general formula is the case of a
one-dimensional Kolmogorov (also known as Fokker-Planck equation)
equation admitting a detailed balance equilibrium $m$. In this case
(\ref{adv-diff1}) takes the form:
$$ \frac{\partial u}{\partial t} -
\frac{\partial}{\partial x} 
\PAR{D(x)m(x)\frac{\partial}{\partial x} \PAR{\frac{u}{m}}}
 = 0.$$
Taking  $v=m$ in (\ref{ent-prod}) we obtain : 
$$
\frac{d}{dt} H^\Phi(u|m) = -\int_{\R^d} \Phi''\PAR{\frac{u}{m}}
\ABS{\frac{\partial}{\partial x} \PAR{\frac{u}{m}}}^2 D m \, dx,
$$
in which one recognizes a generalization of the familiar Fisher
information, the usual Fisher information corresponding to the case
where $D(x)=1$ for all $x$, together with the choice $\Phi(r) = r \log
r$.

Let us now go over to the case of a general linear second-order scalar
advection-diffusion equation:
\begin{equation} \label{adv-diff}
\frac{\partial u}{\partial t} + \hbox{\rm div}
 (b(t,x)u - a(t,x) \nabla u) = 0.
\end{equation}
Here $b$ is a given vector field, and $a$ is a given diffusion matrix.
For any two positive solutions $u,v$ to this equation,
(\ref{ent-prod}) now becomes:
\begin{equation} \label{H-prime}
\frac{d}{d t} H^\Phi(u|v) = 
- \int_{\R^d} \Phi''\PAR{\frac{u}{v}}\SBRA{\nabla \PAR{\frac{u}{v}} \cdot 
                            a \nabla \PAR{\frac{u}{v}}} v \, dx.
\end{equation}
As expected, this formula shows that linear transport does not play
any role in the entropy production (the expression does not involve
the velocity field $b$): in physical parlance, diffusion here is the
only irreversible process.  Let us emphasize that, from the convexity
of $\Phi$ and the positive definiteness of the diffusion matrix $a$,
we obtain that the relative entropy $H(f|g)$ is a time-decreasing
quantity, {\bf whatever the solutions $u,v$ are}, and {\bf whatever
  the coefficients $b,a$ are}. In particular for arbitrary
time-dependent coefficients $b$ and $a$ the system will have no
detailed balance equilibrium , i.e. the problem
$$
b(t,x)u(t,x) = a(t,x) \nabla u(t,x)
$$
will have no solution at all.  Note that in the context of Markov
processes, a similar dissipation property was exhibited by Yosida, and
Kubo (see \cite{K}).

The natural question then arises to investigate under which
(sufficient) conditions on the coefficients of the equation does the
entropy decrease to zero.  In the classical setup where one
investigates the trend toward a stationary solution, it is
well-known that for a large class of such stationary solutions, some
Logarithmic Sobolev Inequality is available \cite{ov}.  This fact can
be used to obtain a Gronwall-type inequality for the entropy, thereby
yielding exponential decay. The supplementary difficulty here is that
the measure relative to which the entropy is computed moves along the
flow, in such a way that classical conditions which ensure that a
logarithmic Sobolev inequality will hold cannot be checked a-priori.
 
In the next section we revisit a well-known prototype, the Ornstein-Uhlenbeck
equation with constant drift.  This example shows that depending on the nature
of the drift the entropy may decay to zero in an exponential or algebraic
fashion, or converge to a nonzero value. Section 3 collects the technical
tools needed to show that for the case where the diffusion matrix is the
identity matrix, the solution of the evolution problem will satisfy the
logarithmic Sobolev inequality at all positive times. The asymptotic behavior
of the entropy is obtained as an easy corollary in Section 4. Finally in
Section 5 we show that at least for the heat equation (but we believe for a
much larger class of parabolic problems), the choice of the fundamental
solution for $v$ provides a very quick proof of the classical Gaussian
intermediate scaling.

\section{The fundamental example: the Ornstein-Uhlenbeck process}

Let us consider the simplest case. Denote by ${(X)}_{t\geq 0}$ the solution of
$$
dX_t=\sqrt{2}B_t-\lambda X_t\,dt,
$$
where ${(B)}_{t\geq 0}$ is a standard Brownian motion on $\R$ and
$\lambda\in\R$ is a constant. This equation can be solved as follows:
$$
X_t=X_0 e^{-\lambda t}+
\sqrt{2}\int_0^t\! e^{\lambda(s-t)}\,dB_s.
$$
As a conclusion, the measure $P_t(\cdot)(x)$ which is defined as the law of
$X_t$ knowing that $X_0=x$ is the Gaussian measure with mean $xe^{-\lambda t}$
and variance $(1-e^{-2\lambda t})/\lambda$. One can then compute the relative
entropy of $P_t(\cdot)(y)$ with respect to $P_t(\cdot)(x)$:
$$
\alpha(t):=H(P_t(\cdot)(y)|P_t(\cdot)(x))
=\frac{\lambda (x-y)^2}{2(e^{\lambda t}-1)},
$$
since $P_t(\cdot)(y)$ and $P_t(\cdot)(x)$ have the same variance.
Of course, in the case when $\lambda=0$, the above formula has to be
understood as
$$
H(P_t(\cdot)(y)|P_t(\cdot)(x))
=\frac{(x-y)^2}{4 t}.
$$

As a conclusion, three different behaviors can occur: 
\begin{itemize}
\item if $\lambda>0$, then $\alpha$ decreases exponentially fast to 0,
  which is natural since $P_t(\cdot)(x)$ converges exponentially fast
  to its invariant measure $\cN(0,1/\lambda)$. 
\item if $\lambda=0$, then $\alpha$ still goes to zero although
  $P_t(\cdot)$ does not converge to a probability measure,
\item if $\lambda<0$, then $\alpha$ converges exponentially fast to a
  nonzero limit: 
  $$
  \alpha(t)=\frac{-\lambda (x-y)^2}{2(1-e^{\lambda t})}
  =-\frac{\lambda}{2}(x-y)^2-\frac{\lambda(x-y)^2}{2(1-e^{\lambda
      t})}e^{\lambda t} \xrightarrow[t\rightarrow
  \infty]{}-\frac{\lambda}{2}(x-y)^2.
$$
\end{itemize}

\section{The local $\Phi$-Sobolev inequality for inhomogeneous
diffusion Semigroups}

\subsection{Notations}

In this section we consider the family of formal elliptic partial
differential operators $(L_t)_{t>0}$ defined by
\begin{equation} \label{def-op}
L_t f(x) := \sum_{i,j=1}^d a_{ij}(t,x)\partial_{ij} f(x) 
+ \sum_{i=1}^d b_i(t,x) \partial_if(x),
\end{equation}
where ${(a_{ij}(t,\cdot))}_{1\leq i,j\leq d}$ is a definite positive diffusion
matrix and $b(t,\cdot)$ is a given vector field on $\R^d$, defined for all
$t>0$. Let us suppose that the coefficients are smooth functions of $(t,x)$.
This family of operators ${(L_t)}_{t\geq 0}$ generates a
inhomogeneous Markov semigroup which we will denote by $(P_{s,t})_{ 0 \leq s \leq t}$ in the
following sense. Writing as usual $a$ as $a=\sigma \sigma^T$, one can
associate to ${(L_t)}_{t\geq 0}$ the solution of the following SDE:
$$
X^{x,r}_t=x+\int_r^t\!b(X^{x,r}_s)\,ds+
\sqrt{2}\int_r^t\!\sigma(X^{x,r}_s)\,dB_t
$$
where ${(B_t)}_{t\geq 0}$ is a standard Brownian motion on $\R^d$.
Semigroup and probabilistic approaches are linked by the fundamental relation
$$
P_{s,t}f(x):=\dE f\PAR{X^{s,x}_t}.
$$
The Markov property of $X$ can be translated into a composition rule for the
semigroup: for every $s\leq t\leq u$, 
$$
P_{s,u}f(x)=\dE\SBRA{f(X^{s,x}_u)}=\dE\SBRA{f\PAR{X^{t,X^{s,x}_t}_u}}
=\dE\SBRA{P_{t,u}f(X^{s,x}_t)}=P_{s,t}P_{t,u}f(x).
$$
This semigroup satisfies the well-known Kolmogorov equations:
\begin{equation} \label{kolmo}
\partial_s P_{s,t} f = -L_s P_{s,t}f, \quad 
\partial_t P_{s,t} f = P_{s,t} L_t f.
\end{equation}

Let $\mu$ be a probability measure on $\dR^d$ and $u$ the density function of
the law of $X_t$ knowing that $\cL(X_0)=\mu$. Then, for every smooth function
$f$,
$$
\dE(f(X_t))=\dE\PAR{\dE(f(X_t)|X_0)}=\int\!P_{0,t}f(x)\,\mu(dx).
$$
The Itô Formula ensures that, for every smooth function $f$ and $s\leq t$,
$$
f(X_t)-f(X_s)-\int_s^t\! L_rf(X_r)\,dr
$$
is a martingale. In other words, 
$$
\dE f(X_t)-\dE f(X_s)-\int_s^t\! \dE L_rf(X_r)\,dr=0.
$$
As a consequence,
$$
\int\!\SBRA{u(t,x)-u(s,x)-\int_s^t\! L_r^*u(r,x)\,dr} f(x)\,dx=0,
$$
and $u$ satisfies Equation \eqref{adv-diff} with initial condition
$\mu$ in a weak sense.

Following \cite{b1,b2}, let us associate to $L_t$ the two bilinear forms
$\Gamma(t)$ and $\gd(t)$ defined by:
\begin{eqnarray*}
\Gamma(t)(f,g) &:=& \frac{1}{2}[ L_t(fg) - g L_t f - f L_t g],\\
\gd(t)(f,g) &:=& \frac{1}{2}[ L_t \Gamma(f,g) - \Gamma(g, L_t f) - \Gamma(f,
L_t g)].
\end{eqnarray*}
We will write $\Gamma(t) (f)$ instead of $\Gamma(t)(f,f)$ and $\gd(t)(f)$
instead of $\gd(t)(f,f)$.

\begin{rmrk}
One can check that
$$
\Gamma(t)(f,g) = \nabla f \cdot a(t,\cdot) \nabla g =
\frac{1}{2} \sum_{i=1}\sum_{j=1} ^d
[a_{ij}+a_{ji}](t,\cdot) \partial_i f\partial_j g.
$$
\end{rmrk}

\begin{rmrk}
  The expression of $\gd$ is much more complicated in the general case. In the
  simple (but informative) case when $a(t,\cdot)$ is the identity matrix, it
  is very easy to check the following formula:
\begin{equation} \label{expression-Gamma}
\gd(t)(f) := ||\text{Hess}(f)||_{2}^2 
- \nabla f \cdot \text{Jac}(b(t)) \nabla f,
\end{equation}
where Hess$(\cdot)$ (resp. Jac$(\cdot)$) stands for the Hessian (resp.
Jacobian) matrix, and $||B||_{2}$ denotes the Hilbert-Schmidt norm.
\end{rmrk}

\begin{rmrk}
  Notice that
  $$
  \nabla f\cdot \text{Jac}(b)\nabla f =\nabla f\cdot \text{SJac}(b)\nabla
  f,
  $$
  where SJac$(\cdot)$ stands for the symmetric part of the Jacobian matrix
  \emph{i.e.},
  $$
  \text{SJac}(b)_{ij}=\frac{\text{Jac}(b)_{ij}+\text{Jac}(b)_{ji}}{2}.
  $$
  The antisymmetric part of the Jacobian of $b$ brings no contribution in
  our study. One can think about the following explicit example: consider the
  2-dimensional process $X$ solution of the following SDE:
  $$
  dX_t=dB_t- \left(
  \begin{array}[c]{cc}
    1 & 1\\
    -1& 1
  \end{array}
  \right) X_t\,dt.
  $$
  The antisymmetric part of the drift induces a rotation whereas the
  symmetric part ensures the convergence to equilibrium.
\end{rmrk}

Let $\Phi$ : $\cI\rightarrow\dR$ be a smooth convex function
defined on a closed interval $\cI$ of $\dR$ not necessarily bounded. Let $\mu$
be a positive measure on a Borel space $(\Omega,\cF)$. The $\Phi$-entropy
functional $\ep{\mu}$ is defined on the set of $\mu$-integrable functions $f$
: $(\Omega,\cF)\rightarrow (\cI,\cB(\cI))$ by
$$
\enp{\mu}{f}=\int_\Omega\!\Phi(f)\,d\mu-\Phi\PAR{\int_\Omega\!f\,d\mu}.
$$
In what follows, $\mu$ is a probability measure. As a
consequence, $\int_\Omega\!f\,d\mu\in \cI$ and the definition make
sense. in the sequel, one has to make an extra assumption in order to
derive interesting functional inequalities:
\begin{equation}
  \label{eq:phisec}
  (u,v)\mapsto \Phi''(u)v^2 \text{ is non negative and convex on } \cI\times\cI. 
\end{equation}

\begin{rmrk}
  The classical variance and entropy are $\Phi$-entropy functionals
  respectively associated to $x\mapsto x^2$ on $\cI=\dR$ and $x\mapsto x\log
  x$ on $\cI=[0,+\infty)$.
\end{rmrk}

\begin{dfntn}
  The semigroup ${(P_{s,t})}_{0\leq s\leq t}$ is said to satisfy a
  local $\Phi$-Sobolev inequality with cons\-tants ${(C_{s,t})}_{0\leq
    s\leq t}$ if for all $s\leq t$ and smooth function $f$,
  $$
  \emph{Ent}^\Phi_{P_{s,t}}(f):=P_{s,t} (\Phi(g)) - \Phi(P_{s,t}g) 
  \leq C_{s,t}P_{s,t}\PAR{\Phi''(f)\Gamma(t)(f)}.
  $$
\end{dfntn}
\begin{rmrk}
  Under the so-called Bakry-Émery criterion, 
  $$
  \exists \rho\in\dR,\quad\forall f\text{ smooth},\quad 
  \gd(f)\geq \rho \Gamma(f),
  $$
  homogeneous diffusion semigroups satisfy a Poincaré and a logarithmic
  Sobolev inequality (see \cite{b3}). As a generalization, $\Phi$-Sobolev
  inequalities can also be established (see \cite{Dj}).
\end{rmrk}

Our aim is to take into account the time dependence of the coefficients of the
diffusion process. We will show that the appropriate adaptation of the
Bakry-Émery criterion to that situation is as follows: 
\begin{equation} \label{curvature}
 \exists \rho\ :\ t\mapsto \rho(t)\in\dR,\quad \forall f\text{ smooth},
 \quad 
 \gd(t)(f)+ \frac{1}{2}\partial_t\Gamma(t)(f)\geq \rho(t) \Gamma(f),
\end{equation}
where  $\partial_t\Gamma(t)$ is defined as
$$
\partial_t\Gamma(t)(f,g):=\sum_{i,j=1}^d\partial_t
a_{ij}(t,\cdot)\partial_i f\partial_j g.
$$

The key point in the homogeneous and diffusive case is to
get the following commutation relation (which turns out to be equivalent to
Bakry-Émery criterion):
$$
\sqrt{\Gamma P_t f}\leq e^{-\rho t}P_t\PAR{\sqrt{\Gamma f}}.
$$
In the following subsection we derive such a commutation relation in the
inhomogeneous case.

\subsection{The commutation relation}

Let $s$ and $t$ be two fixed times, with 
$0\leq s \leq t$.  The key point is the following lemma, which describes
how the dissipative mechanism tends to flatten gradients:

\begin{lmm}\label{lem:comm}
  Suppose that the family of operators ${(L_t)}_{t\geq 0}$ defined in
  (\ref{def-op}) satisfies \eqref{curvature}. For any $\tau$ between 0 and
  $t$, the following inequality holds true:
\begin{equation} \label{flat-gradient}
\sqrt {\Gamma(\tau) (P_{\tau,t}g)} \leq \exp\PAR{{-\int_\tau^t \rho(u) \, du}} 
P_{\tau,t}\PAR{\sqrt {\Gamma(t)(g)}}.
\end{equation}
\end{lmm}

\begin{proofp}{ } 
For all $u \in [\tau,t]$, we define $\beta(u)$ by
$$
\beta(u) = P_{\tau,u}\PAR{ \sqrt{ \Gamma (u)(P_{u,t}g)}},
$$
and compute its derivative by using (\ref{kolmo}). The crucial assumption
that $L_u$ is a diffusion operator ensures that
\begin{equation}
  \label{eq:dif}
L_u(\Phi(g))=\Phi'(g)L_u(g)+\Phi''(g)\Gamma(u)(g).
\end{equation}
In order to make the exposition clearer we denote $P_{u,t}g$ as $h$ and
$\Gamma$ stands for $\Gamma(u)$. A straightforward computation leads to:
\begin{eqnarray*}
\beta'(u) & = & P_{\tau,u}\PAR{ L_u \sqrt{ \Gamma P_{u,t}g}
+ \frac{1}{2 \sqrt{ \Gamma(P_{u,t}g)}}
\BRA{-2\Gamma(P_{u,t}g,L_uP_{u,t}g)
+\partial_u \Gamma(P_{u,t}g)}} \\
& = & P_{\tau,u} \SBRA{\frac{1}{2 \sqrt{ \Gamma h}}L_u \Gamma h 
- \frac{1}{4 (\Gamma h)^{\frac{3}{2}}} \Gamma\Gamma h
- \frac{\Gamma(h,L_uh)}{\sqrt{ \Gamma h}}
+\frac{\partial_u\Gamma(h)}{2\sqrt{ \Gamma h}} } \\
& = & P_{\tau,u}\SBRA{ \frac{ 2 (\Gamma h)  L_u \Gamma h
- 4 (\Gamma h) \Gamma(h,L_uh) - \Gamma\Gamma h +2(\Gamma h)(\partial_u\Gamma
h)} {4 (\Gamma h)^{\frac{3}{2}}}} \\
& = & P_{\tau,u} \SBRA{\frac{ 4 (\Gamma h) \gd(u)(h) - \Gamma\Gamma h+2(\Gamma
  h)\partial_u\Gamma h}
{4 (\Gamma h)^{\frac{3}{2}}}}.
\end{eqnarray*}

Therefore we obtain:
$$
\beta'(u) - \rho(u) \beta(u) = P_{\tau,u} \SBRA{\frac {4(\Gamma h)(\gd(u)(h)
    +(1/2)\partial_u\Gamma h- \rho(u)\Gamma h) - \Gamma\Gamma h} {4 (\Gamma
    h)^{\frac{3}{2}}}}.
$$
Following \cite{b3}, one can show, thanks to the diffusion assumption, that
the criterion \eqref{curvature} implies that, for all smooth functions $f$, 
$$
\gd(u)(f)+(1/2)\partial_u\Gamma f- \rho(u)\Gamma f \geq \frac{\Gamma\Gamma
  f}{4\Gamma f}.
$$
\begin{rmrk}
  In the case when $a$ is the identity matrix, one can easily get derive this
  inequality from the criterion \eqref{curvature} and the Cauchy-Schwarz
  inequality.
\end{rmrk}

As a conclusion, $\beta$ satisfy the following differential inequality:
$$\beta'(u) \geq  \rho(u) \beta(u).$$
In other words, the function 
$$
u\mapsto \beta(u)\exp\PAR{-\int_\tau^u\!\rho(v)\,dv}
$$
is an increasing function on the interval $[\tau,t]$ which implies
that
$$
\beta(\tau) \leq \beta(t) \exp\PAR{-\int_\tau^t \rho(u) \, du}.
$$
This is precisely the desired inequality.
\end{proofp}

\subsection{Local $\Phi$-Sobolev inequalities}

\begin{thrm}
  Suppose that the family of operators ${(L_t)}_{t\geq 0}$ defined in
  (\ref{def-op}) satisfies \eqref{curvature}. Then for any times $s,t$
  with $0 \leq s \leq t$ and any positive function $g$, $P_{s,t}$
  satisfies the following $\Phi$-Sobolev inequality:
  $$
  \emph{Ent}^\Phi_{P_{s,t}}(g)\leq c(s,t)
  P_{s,t}\PAR{\Phi''(g)\Gamma(t)(g)},
  $$
  where the constant $c(s,t)$ can be chosen as:
  $$
  c(s,t) = \int_s^t \exp\PAR{-2 \int_\tau^t \rho(u) \, du} \, d \tau.
  $$
\end{thrm}

\begin{proof}
Consider the function $\alpha$ : $[s,t]\rightarrow \dR$ defined by :
$$
\alpha(\tau):= P_{s,\tau} \PAR{\Phi(P_{\tau,t}g)},
$$
Let us compute the derivative of $\alpha$: 
$$
\alpha'(\tau)= P_{s,\tau} \PAR{ L_\tau \PAR{\Phi(P_{\tau,t}g)}
-\Phi'(P_{\tau,t}g)L_\tau P_{\tau,t}g}
$$
Thanks to the fact the diffusion assumption, \eqref{eq:dif} ensures
that 
$$
\alpha'(\tau)= P_{s,\tau}
\PAR{\Phi''(P_{\tau,t}g)\Gamma(\tau)(P_{\tau,t}g)}.
$$
The commutation relation \eqref{flat-gradient} ensures that 
\begin{eqnarray*}
\Phi''(P_{\tau,t}g)\Gamma(\tau)(P_{\tau,t}g)
&=&\Phi''(P_{\tau,t}g)\PAR{\sqrt{\Gamma(\tau)(P_{\tau,t}g)}}^2\\
&\leq &\exp\PAR{-2\int_\tau^t\!\rho(u)\,du}
\Phi''(P_{\tau,t}g) P_{\tau,t}\PAR{\sqrt{\Gamma(t)(g)}}^2.
\end{eqnarray*}
Jensen inequality with the bivariate function
$(u,v)\mapsto\Phi''(u)v^2$ (which is assumed to be convex according to
\eqref{eq:phisec}) ensures that
$$
\Phi''(P_{\tau,t}g) P_{\tau,t}\PAR{\sqrt{\Gamma(t)(g)}}^2
\leq  P_{\tau,t}\PAR{\Phi''(g)\Gamma(t)(g)}.
$$
As a conclusion, 
$$
\alpha'(\tau)\leq \exp\PAR{-2\int_\tau^t\!\rho(u)\,du}
P_{s,\tau}P_{\tau,t}\PAR{\Phi''(g)\Gamma(t)(g)}
=\exp\PAR{-2\int_\tau^t\!\rho(u)\,du}
P_{s,t}\PAR{\Phi''(g)\Gamma(t)(g)} .
$$
Since $\alpha(t)=P_{s,t}(\Phi(g))$ and
$\alpha(s)=\Phi(P_{s,t}(g))$, the result follows upon integration of
this inequality between $s$ and $t$.
\end{proof}

In the special case of $\Phi$ : $x\mapsto x\log x$, the
local logarithmic Sobolev inequality for the semigroup ${(P_{s,t})}_{0\leq
  s\leq t}$ can be stated as follows:

\begin{crllr}\label{ze-result}
  Suppose that the family of operators ${(L_t)}_{t\geq 0}$ defined in
  (\ref{def-op}) satisfies \eqref{curvature}. Then for any times $s,t$ with $0
  \leq s \leq t$ and any positive function $g$, $P_{s,t}$ satisfies the
  following logarithmic Sobolev inequality:
\begin{equation} \label{LISL}
\emph{Ent}_{P_{s,t}}(g):=P_{s,t} (g \log g) - (P_{s,t}g) \log (P_{s,t}g) 
\leq c(s,t) P_{s,t}\PAR{\frac{\Gamma(t)(g)}{g}},
\end{equation}
where the constant $c(s,t)$ can be chosen as:
$$
c(s,t) = \int_s^t \exp\PAR{-2 \int_\tau^t \rho(u) \, du} \, d \tau.
$$
\end{crllr}

\begin{rmrk}
If for every $x\in \dR^d$, the matrix ${(a_{ij}(t,x))_{i,j}}$ is bounded by
the identity matrix (as symmetric bilinear forms), then
$\Gamma(t)(g)\leq\ABS{\nabla g}^2$ and $P_{s,t}$ satisfies the
classical logarithmic Sobolev inequality: 
$$
\emph{Ent}_{P_{s,t}}(g)
\leq c(s,t) P_{s,t}\PAR{\frac{\ABS{\nabla g}^2}{g}}.
$$
\end{rmrk}

\begin{rmrk}
  In the case where $\rho(t)=\rho$ for all $t>0$, we recover the 
  constant 
  $$
  c(s,t) = \frac{1}{2 \rho} \PAR{1 - e^{-2 \rho(t-s)}}
  $$
  provided by \cite{b3} in the homogeneous case.
\end{rmrk}

\section{The Kullback-Leibler distance of two orbits of a parabolic problem}

We now consider the following parabolic equation \eqref{adv-diff}:
\begin{eqnarray}
\frac{\partial v}{\partial t}(t,x) &=& L_t^*v(t,x), 
\quad t>0, \, x \in \R^d, \label{ze-pde} \\
v(0,x) &=& v_0(x), \quad x \in \R^d,\label{ze-id}
\end{eqnarray}
where $L_t^*$ stands for the adjoint of $L_t$ with respect to the Lebesgue
measure on $\dR^d$. 

If we assume that the initial data $v_0$ satisfies a logarithmic
Sobolev inequality, the result of the previous section may be used to
show that the inequality is propagated in time:

\begin{thrm} \label{propag}
  Assume the initial data $v_0$ satisfies the following logarithmic Sobolev
  inequality: for all smooth functions $f$, 
$$
\int\! f\log f v_0-\int\! fv_0\log\int\! fv_0
\leq d_0\int\!\Gamma(0)(f)v_0. 
$$
If the family ${(L_t)}_{t\geq 0}$ satisfies the criterion
\eqref{curvature}. Then for any positive time, the solution $v(t,\cdot)$ of
(\ref{ze-pde})(\ref{ze-id}) satisfies the following logarithmic Sobolev
inequality 
$$
\int\! f\log f v(t,\cdot)-\int\! fv(t,\cdot)\log\int\! fv(t,\cdot)
\leq d(t)\int\!\Gamma(t)(f)v(t,\cdot),
$$
where
\begin{equation} \label{la-constante}
d(t) := d_0 \exp\PAR{-2 \int_0^t \rho(r) \, dr} + 
\int_0^t \exp\PAR{-2 \int_\tau^t \rho(r) \, dr} \, d \tau.
\end{equation}
\end{thrm}

\begin{proofp}{}
For any positive function $g$ we have:
$$
\int g(x) \log g(x) v(t,x) \, dx 
=\int\!\dE\SBRA{f\PAR{X^{0,x}_t}}\,v_0(x)\,dx
= \int P_{0,t}(g \log g)(x) v(0,x) \, dx. 
$$
The integrand may be bounded from above by using the local inequality
\eqref{LISL} with $s=0$, and the logarithmic Sobolev inequality for $v_0$:
\begin{eqnarray*}
\int P_{0,t}(g \log g) v(x,0) \, dx & \leq & 
\int P_{0,t} g \log P_{0,t} g v_0(x) \, dx 
+ c(0,t)\int P_{0,t}\PAR{\frac{\Gamma(t)(g)}{g}} v_0(x) \, dx  \\ 
& \leq & d_0\int  \frac{\Gamma(0)(P_{0,t}g)}{P_{0,t} g} v_0(x) \, dx
+ c(0,t) \int P_{0,t}\PAR{\frac{\Gamma(t)(g)}{g}} v_0(x) \, dx .
\end{eqnarray*}
The first integral in the last inequality may be estimated by
\eqref{flat-gradient}, and this completes the proof.
\end{proofp}

We are now in a position to estimate the Kullback-Leibler distance
between two arbitrary orbits of (\ref{ze-pde}):
\begin{thrm}
  Under the same assumptions as in previous theorem, let $u$ be
  another solution of (\ref{ze-pde}) (i.e., corresponding to different
  initial data $u_0$).  Assume that $u_0$ and $v_0$ are positive.
  Define the relative entropy of $u$ with respect to $v$ at any
  positive time $t$ by
$$H(u(t)|v(t)) := \int u(t,x) \log \frac{u(t,x)}{v(t,x)} \, dx.$$
This quantity is then bounded as follows for all positive times:
$$H(u(t)|v(t)) \leq H(u(0)|v(0)) c(t),$$
where 
\begin{equation} \label{def-c}
c(t)=\exp\PAR{-\int_0^t \frac{1}{d(s)} \, ds},
\end{equation}
and $d(t)$ is the constant defined in (\ref{la-constante}).
\end{thrm}

\begin{proofp}{}
The proof is a straightforward application of Gronwall's lemma.
Let us set $g:=\frac{u}{v}$; from (\ref{H-prime}) with $\Phi(z) = z
\log z$, we obtain
$$
\frac{d}{dt}H(u(t)|v(t)) = - \int \frac{\Gamma(t)(g)}{g} v \, dx.
$$
Therefore, Theorem \ref{propag} gives the following control:
$$
\frac{d}{dt}H(u(t)|v(t)) \leq -\frac{1}{d(t)} H(u(t)|v(t)),
$$
which gives the result.  Let us conclude this section by indicating
what the obtained rate is when the quantity $\rho(t)$ may be taken to
be a fixed constant $\rho$.  If $\rho=0$, then (\ref{la-constante})
gives $d(t)=d_0+t$, therefore we get the algebraic decay
$$
c(t) = \frac{1}{1+\frac{t}{d_0}}.
$$
For $\rho \neq 0$, we have 
$$d(t)= d_0 e^{-2 \rho t} + \frac{1}{2\rho}(1 - e^{-2 \rho t}),
$$
and the integral in (\ref{def-c}) may be computed to yield
$$
c(t) = 
\frac{2 \rho d_0  e^{-2 \rho t}}{1 + (2 \rho d_0 -1)e^{-2 \rho t}}.
$$
\end{proofp}

\section{An application to intermediate asymptotics}

In the case where a Kolmogorov equation has a rather trivial (e.g.
constant) asymptotic state, it is often the case that (due to the self
similarity of the underlying Markov process) some appropriate
rescaling of the orbit shows structure (e.g. Gaussian), a phenomenon
termed {\it intermediate asymptotics} by Barenblatt \cite{B}. To
illustrate this point let us consider the one dimensional linear heat
equation on the entire line:
$$\frac{\partial u}{\partial t} = \frac{1}{2} \frac{\partial^2 u}{
  \partial t^2}, \quad x \in \R.
$$
It is very easy to check that for large times $u$ converges
point wise to zero, but it is also very easy to read off from the
explicit form of the solution (as given in terms of the heat kernel)
the following point wise convergence:
\begin{equation}
\sqrt t u(y \sqrt t ,t) \rightarrow C \exp\PAR{- \frac{y^2}{2}},
\label{int-as}
\end{equation}
where the constant $C$ may be determined from mass conservation.  As
is well-known, this convergence may be obtained by rescaling the
equation and constructing an entropy functional for the rescaled
equation.  Let us briefly recall this argument. The first step
consists in rescaling the function $u$ by setting:
$$
u(t,x) = \alpha(t)v\PAR{ \frac{x}{\beta(t)},\tau(t)},
$$
where the scaling functions $\alpha(t),\beta(t),\tau(t)$ are to be
chosen so as to make the equation for $v$ as simple as possible, while
preserving the mass constraint:
$$
\int_{\R^d} u(t,x) \, dx = const.
$$
These two requirements lead to the choice
\begin{equation}
v(y,\ln t) = \sqrt t u(y \sqrt t,t),
\label{rescaled}
\end{equation}
thus to the equation
$$
\frac{\partial v}{\partial \tau} = \frac{1}{2} \frac{\partial}{\partial y}
\PAR{yv + \frac{\partial v}{\partial y}}.
$$
The detailed balance equilibria for this equation are exactly all
multiples of the standard Gaussian density, i.e. take the form:
\begin{equation} \label{fixed-point}
v_\infty(y) = C \exp\PAR{- \frac{y^2}{2}},
\end{equation}
for some constant $C$.  By the result of previous section, the
relative entropy
\begin{equation}
H(v|v_\infty) = \int_\R v \ln \frac{v}{v_\infty} \, dy
\label{resc-ent}
\end{equation}
decreases to zero for large times, which implies that $v$ converges to
a fixed point of the form (\ref{fixed-point}) for some $C$.

This result, when rephrased in terms of $u$, is exactly the
intermediate asymptotics (\ref{int-as}).  Note that in this case the
constant $C$ is uniquely determined from the mass conservation
relation:
$$
\int_\R u(t,x) \, dx = \int_\R u(x,0) \, dx.
$$
The key point is now that by using the change of variable $x = y
\sqrt t$ in (\ref{resc-ent}) one obtains:
$$
H(v|v_\infty) = \int_\R u(t,x) \ln \frac{u(t,x)}{C \frac{1}{\sqrt t} 
\exp{(- \frac{x^2}{2t})}} \, dx,
$$
which is exactly the relative entropy of $u$ with respect to the
fundamental solution of the heat equation. The fact that this entropy
is dissipated combined with Pinsker's inequality now immediately leads
to the intermediate asymptotic (\ref{int-as}), without having to
resort to any rescaling of the function $u$.

In other words, the fundamental solution encodes the intermediate
asymptotics.

\end{document}